\documentclass[12pt]{amsart}
\usepackage{amsfonts}
\usepackage{amssymb}
\usepackage{amsmath}
\usepackage{amsthm}

\normalsize

\newtheorem{theorem}{Theorem}

\newtheorem{proposition}{Proposition}

\newtheorem{corollary}{Corollary}

\theoremstyle{definition}
\newtheorem{definition}{Definition}

\theoremstyle{definition}

\begin{document}

\title[Function Models for Teichm\"uller Spaces]{Function Models for Teichm\"uller Spaces
and Dual Geometric Gibbs Type Measure Theory for Circle Dynamics}

\author{Yunping Jiang}

\address{Department of Mathematics\\
Queens College of the City University of New York\\
Flushing, NY 11367-1597\\
and\\
Department of Mathematics\\
Graduate School of the City University of New York\\
365 Fifth Avenue, New York, NY 10016}
\email{yunping.jiang@qc.cuny.edu}



\thanks{The research is partially supported by
PSC-CUNY awards.}
\thanks{This article is prepared for the proceedings of the International Workshop on
Teichmueller theory and moduli problems held at the Harish-Chandra
Research Institute (HRI), Allahabad, India, January 5 to 14, 2006.}

\begin{abstract}
Geometric models and Teichm\"uller structures have been introduced
for the space of smooth expanding circle endomorphisms and for the
space of uniformly symmetric circle endomorphisms. The latter one
is the completion of the previous one under the Techm\"uller
metric. Moreover, the spaces of geometric models as well as the
Teichm\"uller spaces can be described as the space of H\"older
continuous scaling functions and the space of continuous scaling
functions on the dual symbolic space. The characterizations of
these scaling functions have been also investigated. The Gibbs
measure theory and the dual Gibbs measure theory for smooth
expanding circle dynamics have been viewed from the geometric
point of view. However, for uniformly symmetric circle dynamics,
an appropriate Gibbs measure theory is unavailable, but a dual
Gibbs type measure theory has been developed for the uniformly
symmetric case. This development extends the dual Gibbs measure
theory for the smooth case from the geometric point of view. In
this survey article, We give a review of these developments which
combines ideas and techniques from dynamical systems,
quasiconformal mapping theory, and Teichm\"uller theory. There is
a measure-theoretical version which is called $g$-measure theory
and which corresponds to the dual geometric Gibbs type measure
theory. We briefly review it too.
\end{abstract}

\maketitle

\newpage

\vspace*{1in}

\centerline{\bf \Large Contents}

\vspace*{0.5in}

\noindent {1.  Introduction.}

\noindent {2.  Circle endomorphisms.}

\noindent {3. Topological models.}

\noindent {4. Geometric models, part I.}

\noindent {5.  Teichm\"uller structures, part I.}

\noindent {6. Characterizations, part I.}

\noindent {7. Teichm\"uller structures, part II.}

\noindent {8. Geometric models, part II.}

\noindent {9. Characterizations, part II.}

\noindent {10. Invariant measures and dual invariant measures.}

\noindent {11. $g$-measures.}

\noindent {12. Geometric Gibbs measures and dual geometric Gibbs
measures.}

\noindent {13. Dual geometric Gibbs type measures.}

\newpage

\subsection*{1. Introduction}$\\$

\vspace*{-5pt}

Circle maps are basic elements in dynamical systems. The dynamics
of a smooth expanding circle map presents many profound phenomena
in mathematics and physics such as the structural stability
theory, ergodic theory, probability theory, and, more recently,
chaos theory. Teichm\"uller theory studies complex manifold
structures of almost complex structures on Riemann surfaces. We
have brought in some concepts and techniques in Teichm\"uller
theory into the study of geometric structures of spaces of circle
expanding maps. In this survey article we review some development
in this direction.

\subsection*{2. Circle endomorphisms}$\\$

\vspace*{-5pt}

The theme here is an orientation-preserving covering map $f$ from
the unit circle $T=\{ z\in {\Bbb C}\; |\; |z|=1\}$ onto itself. Let
$d$ be the topological degree of $f$. We assume that $d\geq 2$. The
universal cover of $T$ is the real line ${\Bbb R}$ with a covering
map
$$
\pi (x) = e^{2\pi x}: {\Bbb R} \to T
$$
Then $f$ can be lifted to an orientation-preserving homeomorphism
$F$ of ${\Bbb R}$ with the property that $F(x+1)=F(x)+d$. Since a
covering map of degree $\geq 2$ has a fixed point, we assume that
$z=1$ is a fixed point of $f$. Then by assuming $F(0)=0$, we set up
a one-to-one correspondence between degree $d$ circle covering maps
$f$ with $f(1)=1$ and real line homeomorphisms $F$ with
$F(x+1)=F(x)+d$. Thus we call in this paper $f$ or the corresponding
$F$ a circle endomorphism. We use $f^{n}$ (or $F^{n}$) to mean the
composition of $f$ (or $F$) by itself $n>0$ times.

A circle endomorphism $f$ is $C^{k}$ for $k\geq 1$ if its
$k^{th}$-derivative $F^{(k)}$ is continuous and $C^{k+\alpha}$ for
some $0< \alpha \leq 1$ if, furthermore, $F^{(k)}$ is
$\alpha$-H\"older continuous, that is,
$$
\sup_{x\neq y\in {\Bbb R}} \frac{|F^{(k)} (x) -F^{(k)}
(y)|}{|x-y|^{\alpha}} =\sup_{x\neq y\in [0,1]} \frac{|F^{(k)} (x)
-F^{(k)} (y)|}{|x-y|^{\alpha}} <\infty.
$$
A $C^{1}$ circle endomorphism $f$ is called expanding if there are
constants $C>0$ and $\lambda >1$ such that
$$
(F^{n})'(x) \geq C\lambda^{n}, \quad n=1, 2, \cdots.
$$

\vspace*{10pt}

\subsection*{3. Topological models}$\\$

\vspace*{-5pt}

The topological classification of smooth expanding circle
endomorphisms was first considered by Shub~\cite{Sh} in 1960's. He
proved that two $C^{2}$ expanding circle endomorphisms $f$ and $g$
are topologically conjugate if and only if they have the same
degree. Here $f$ and $g$ are topologically conjugate if there is a
homeomorphism $h$ of $T$ such that
$$
f\circ h =h\circ g.
$$
By also considering the lift $G$ of $g$ and the lift $H$ of $h$, we
have an equivalent definition that $f$ and $g$ are topologically
conjugate if there is a homeomorphism $H$ of ${\Bbb R}$ with
$H(x+1)=H(x) +1$ such that
$$
F\circ H = H\circ G \pmod{1}.
$$

Shub's proof is an application of the contracting fixed point
theorem in functional analysis. Consider the space ${\mathcal C}$ of
all continuous function $\phi$ on ${\Bbb R}$ with $\phi(x+1)
=\phi(x) +1$ with the maximum norm
$$
||\phi|| = \sup_{x\in [0,1]} |\phi(x)|.
$$
Then ${\mathcal C}$ is a Banach space. Define an operator ${\mathcal
L} ={\mathcal L}_{F, G}$ as
$$
{\mathcal L}\phi (x) = F^{-1} \circ \phi \circ G : {\mathcal C}\to
{\mathcal C}.
$$
(Note that $F^{-1} (x+d) =F^{-1}(x)+1$.)  Without loss of
generality, we assume that $C=1$. Then one can check that
$$
||{\mathcal L}\phi -{\mathcal L}\psi|| \leq \frac{1}{\lambda} ||
\phi-\psi||.
$$
So ${\mathcal L}$ is a contracting functional from the Banach space
${\mathcal C}$ into itself. And thus it has a unique fixed point
$H$, that is,
$$
F^{-1} \circ H \circ G  =H
$$
Similarly, ${\mathcal L}_{G,F}$ has a unique fixed point
$\widetilde{H}$, that is,
$$
G^{-1} \circ \widetilde{H} \circ F  =\widetilde{H}.
$$
This implies that
$$
H\circ \widetilde{H} =id.
$$
Therefore, $H$ is a homeomorphism of ${\mathcal R}$ such that
$$
F\circ H = H\circ G \pmod{1}.
$$

There is a more general theorem if we bring in the consideration of
Markov partitions. Consider a partition of $[0,1]$ by
$$
I_{i}=I_{i,f}=F^{-1}([i, i+1]), \quad 0\leq i\leq d-1.
$$
It is a Markov partition in the following sense: If we consider a
corresponding partition of $T$, which we still denote as $\{
I_{i}\}_{i=0}^{d-1}$, then
\begin{enumerate}
\item the union of these intervals is $T$;
\item all the intervals in the partition have pairwise disjoint interiors;
\item the restriction of $f$ to the interior of every interval in
the partition is injective.
\end{enumerate}
We use
$$
\eta_{0}=\eta_{0,f}=\{ I_{i}\}_{i=0}^{d-1}
$$
to denote this initial Markov partition. We then have a sequence of
Markov partitions
$$
\eta_{n} =\eta_{n,f}= f^{-n}(\eta_{0}), \quad n=0,1, 2,\cdots
$$
on the unit circle $T$ as well as the unit interval $[0,1]$. We can
label each interval in $\eta_{n}$ as follows. Define
$$
g_{i} (x) = F^{-1} (x+i): [0,1]\to I_{i}, \quad i=0, 1, \cdots,
d-1.
$$
Each $g_{i}$ is a homeomorphism. Given a word $w_{n} =i_{0}\cdots
i_{k}\cdots i_{n-1}$ of $\{ 0, \cdots, d-1\}$ of length $n\geq 1$,
define
$$
g_{w_{n}} =g_{i_{0}}\circ g_{i_{1}}\circ \cdots \circ g_{i_{n-1}}.
$$
Let
$$
I_{w_{n}} = I_{w_{n},f}= g_{w_{n}}([0,1]).
$$
Then
$$
\eta_{n} =\{ I_{w_{n}}\; |\; w_{n}=i_{0}\cdots i_{k}\cdots i_{n-1},
\; i_{k}\in \{ 0, \cdots, d-1\}\}.
$$
One can check that for a word $w=i_{0}\cdots i_{n-1}i_{n} \cdots$ of
infinite length, and with $w_{n} =i_{0}\cdots i_{n-1}$, then
$$
\cdots \subset I_{w_{n}} \subset I_{w_{n-1}} \subset \cdots
I_{w_{1}}\subset [0,1].
$$
Since each $I_{w_{n}}$ is compact,
$$
I_{w}=\cap_{n=1}^{\infty} I_{w_{n}} \neq \emptyset.
$$
Consider the space
$$
\Sigma^{+}=\Sigma^{+}_{d} =\prod_{n=0}^{\infty} \{ 0, 1, \cdots,
d-1\}
$$
$$
=\{ w=i_{0}i_{1}\cdots i_{k} \cdots i_{n-1} \cdots \; |\; i_{k}\in
\{ 0, 1, \cdots, d-1\}, \; k=0, 1, \cdots \}
$$
with the product topology. Then it is a compact topological space.

If each $I_{w}=\{ x_{w}\}$ contains only one point, then we define
the projection $\pi_{+}=\pi_{+,f}$ from $\Sigma^{+}$ onto $T$ as
$$
\pi_{+} (w) =x_{w}.
$$
The projection $\pi_{+}$ is 1-1 except for a countable set $B$
consisting of all labellings $w$ of endpoints in the partitions
$\eta_{n} =\{ I_{w_{n}}\}$, $n=0, 1,\cdots$.

Let
$$
\sigma^{+} (w) = i_{1}\cdots i_{n-1}i_{n} \cdots
$$
be the left shift map. Then $(\Sigma^{+}, \sigma^{+})$ is called a
symbolic dynamical system. From our construction, one can check that
$$
\pi_{+}\circ \sigma^{+} (w) = f \circ \pi_{+} (w), \quad w\in
\Sigma^{+}.
$$
Let
$$
\varepsilon_{n}= \varepsilon_{n,f}=\max_{w_{n}} |I_{w}|
$$
where $w_{n}$ runs over all words of length $n$ of $\{ 0, 1, \cdots,
d-1\}$. Then we have a more general Shub type theorem.

\vspace*{10pt}
\begin{theorem}~\label{tc} Let $f$ and $g$ be two circle endomorphisms such that
both $\varepsilon_{n,f}$ and $\varepsilon_{n,g}$ tend to zero as
$n\to \infty$. Then $f$ and $g$ are topologically conjugate if and
only if their topological degrees are the same.
\end{theorem}

\begin{proof} Since both sets $I_{w,f}=\{x_{w}\}$ and
$I_{w,g}=\{y_{w}\}$ contain only a single point for each $w$, we
define
$$
h (x_{w}) =y_{w}.
$$
One can check that $h$ is a homeomorphism with the inverse
$h^{-1}(y_{w})=x_{w}$.
\end{proof}

Therefore, for a fixed degree $d>1$, there is only one topological
model $(\Sigma^{+}, \sigma^{+})$  for the dynamics of all circle
endomorphisms of degree $d$ with $\varepsilon_{n}\to 0$.

\vspace*{10pt}

\subsection*{4. Geometric models, part I}$\\$

\vspace*{-5pt}

The next theme is the study of geometric models. A result
analogous to Mostow's rigidity theorem for closed hyperbolic
$3$-manifolds was proved by Shub and Sullivan~\cite{SS}. The
result can be stated as follows: Suppose $f$ and $g$ are two
topologically conjugate real analytic expanding circle
endomorphisms. If the conjugacy $h$ is absolutely continuous, it
must be also real analytic. Later, this result was proved for a
more general case: Suppose $f$ and $g$ are two topologically
conjugate $C^{k+\alpha}$ expanding circle endomorphisms for $1\leq
k\leq \omega$ and $0<\alpha \leq 1$. If the conjugacy $h$ is
absolutely continuous, it must be also $C^{k+\alpha}$. Smooth
invariants of a circle endomorphism have also been investigated. A
quantity is called a smooth invariant if it is the same for $f$
and $g$ as long as $f$ and $g$ are smoothly conjugate (this means
that the conjugacy is $C^{k}$ for $k\geq 1$). A point $p$ of $f$
is called a periodic point of period $n\geq 1$ if $f^{i}(p) \neq
p$ for $0\leq i\leq n-1$ but $f^{n}(p)=p$. The eigenvalue at a
periodic point $p$ of period $n$ is defined as $e_{p}
=(f^{n})'(p)$. The eigenvalue $e_{p}$ is a smooth invariant. The
set of all eigenvalues of a $C^{1+\alpha}$ expanding circle
endomorphism is actually a set of complete smooth invariants,
where $0<\alpha \leq 1$. This means that two $C^{1+\alpha}$
expanding circle endomorphisms $f$ and $g$ of degree $d>1$ are
smoothly conjugate if and only if their eigenvalues at the
corresponding periodic points are the same. Therefore, one can use
the set of all eigenvalues to classify geometric models of smooth
expanding circle endomorphisms of the same degree. (Research in
this direction has been extended to a larger class which even
allows one to include maps with critical points. The reader who is
interested in the smooth classification of one-dimensional
dynamical systems in this direction may refer
to~\cite{Ji1,Ji2,Ji3,Ji4,Ji5} for more details.)

However, the structure of the set of all eigenvalues is not clear.
In what follows, we define a function which is called a scaling
function and will contain full information about the set of all
eigenvalues in this context. (The name of the scaling function in
this context was first used by Feigenbaum~\cite{Fe} in describing
the universal geometric structure of attractors of infinitely
period doubling folding maps. It was then used by
Sullivan~\cite{Su} for Cantor sets on the line to describe
differential structures for fractal sets. The present form of the
definition was formulated in~\cite{Ji1} for any Markov map and
then used to study the smooth classification of one-dimensional
maps which have certain Markov properties. The reader may also
refer to~\cite{Ji2,Ji3} for more details.)

As we have already seen, given a circle endomorphism of degree
$d>1$, there is an interval system
$$
\{ \eta_{n}\}_{n=1}^{\infty}  =\{ \{ I_{w_{n}}\}_{w_{n}}
\}_{n=1}^{\infty}
$$
where $w_{n}$ runs over all words of length $n$ of $\{ 0, 1,
\cdots, d-1\}$. When we constructed the topological model
$(\Sigma^{+}, \sigma^{+}\}$ from this interval system, we read
each $w_{n}$ from the left to the right, i.e.,
$w_{n}=i_{0}i_{1}\cdots i_{n-1}$. From the topological point of
view, this means that we consider the set of all left cylinders
$$
[w_{n}]= [w_{n}]_{l} = [i_{0} i_{1}\cdots i_{n-1}]_{l} = \{
w'=i_{0}'i_{1}'\cdots i_{n-1}'i_{n}' \cdots \; |\; i_{0}'=i_{0},
\cdots, i_{n-1}'=i_{n-1}\}
$$
as a basis for the topology.

Now let us consider another topology which has a basis consisting
of all right cylinders. We read $w_{n}$ from the right to the
left, $\kappa_{n}=w_{n}=j_{n-1}j_{n-2}\cdots j_{0}$ and define
$$
\Sigma^{-} = \Sigma^{-}_{d}=\{\kappa =\cdots j_{n-1}\cdots
j_{k}\cdots j_{1}j_{0} \; |\; j_{k}\in \{ 0, 1, \cdots, d-1\}, \;
k=0, 1, \cdots \}.
$$
It is a topological space with a basis for the topology consisting
of all right cylinders
$$
[\kappa_{n}] =[\kappa_{n}]_{r} =[j_{n-1}\cdots j_{0}]_{r} = \{
\kappa'=\cdots j_{n}'j_{n-1}'\cdots j_{0}' \; |\; j_{n-1}'=j_{n-1},
\cdots, j_{0}'=j_{0}\}.
$$
Consider the right shift map
$$
\sigma^{-} : \cdots j_{n-1}\cdots j_{1}j_{0} \mapsto  \cdots
j_{n-1}\cdots j_{1}.
$$
Then we call $(\Sigma^{-}, \sigma^{-})$ the dual symbolic dynamical
system for $f$.

Another way to view the symbolic dynamical system and the dual
symbolic dynamical system is to consider the inverse limit of $f:
T\to T$. This inverse limit can be viewed as a solenoid with the
symbolic representation
$$
\Sigma =\Sigma^{-}\times \Sigma^{+}.
$$
Then $\Sigma^{-}$ represents the transversal direction and
$\Sigma^{+}$ represents the leaf direction.

On the transversal direction $\Sigma^{-}$, we define a function
called the scaling function for $f$ as follows. For any
$\kappa=\cdots j_{n-1}\cdots j_{1}j_{0}\in \Sigma^{-}$, let
$\kappa_{n} =j_{n-1}\cdots j_{1}j_{0}$, $\sigma^{-} (\kappa_{n}) =
j_{n-1}\cdots j_{1}$. Then
$$
I_{\kappa_{n}} \subset I_{\sigma^{-} (\kappa_{n})}.
$$
Define
$$
S(\kappa_{n}) =S_{f}
(\kappa_{n})=\frac{|I_{\kappa_{n}}|}{|I_{\sigma^{-} (\kappa_{n})}|}.
$$

\vspace*{10pt}
\begin{definition}
If for every $\kappa\in \Sigma^{-}$,
$$
S(\kappa) = S_{f}(\kappa)= \lim_{n\to \infty} S(\kappa_{n})
$$
exists, then we have a function
$$
S=S_{f}: \Sigma^{-}\to {\Bbb R}^{+}.
$$
We call this function the scaling function of $f$.
\end{definition}

The space $\Sigma^{-}$ is a metric space with the metric
$$
d(w, w') =\sum_{k=0}^{\infty} \frac{|i_{k}-i_{k}'|}{d^{k}}.
$$
A function $S$ on $\Sigma^{-}$ is called H\"older continuous if
there are constants $C>0$ and $0<\beta\leq 1$ such that
$$
|S(\kappa)-S(\kappa')| \leq C\big( d(\kappa,\kappa')\big)^{\beta},
\quad \kappa, \kappa'\in \Sigma^{-}.
$$

Let ${\mathcal C}^{1+}$ denote the space of all $C^{1+\alpha}$
expanding circle endomorphisms for some $0<\alpha\leq 1$. We have:

\vspace*{10pt}
\begin{theorem}~\label{sc}
The scaling function $S$ for $f\in {\mathcal C}^{1+}$ exists and is
a H\"older continuous function. Furthermore, $S$ is a completely
smooth invariant. This means that $f, g\in {\mathcal C}^{1+}$ are
$C^{1}$ conjugate if and only if they have the same scaling
functions, i.e., $S_{f}=S_{g}$.
\end{theorem}

This result is actually proved for a larger class of
one-dimensional maps which may have critical points. The reader
who is interested in this direction can refer
to~\cite{Ji2,Ji3,Ji4,Ji5}.

Thus geometric models of ${\mathcal C}^{1+}$ can be represented by
degrees $d>1$ and scaling functions $S$ as follows. We say
$f\sim_{s} g$ if $f$ and $g$ are $C^{1}$ conjugate. It is an
equivalence relation in ${\mathcal C}^{1+}$. Then we have that
$$
{\mathcal C}^{1+} /\sim_{s}  =\{d, S\}.
$$
For a fixed $d>1$, let ${\mathcal C}^{1+}_{d}$ be the space of $f\in
{\mathcal C}^{1+}$ with the degree $d$. Then
$$
{\mathcal C}^{1+}_{d} /\sim_{s}  =\{S\}.
$$

There are two natural problems now. One is to study the geometric
structure on ${\mathcal C}^{1+}_{d} /\sim_{s}$. The other is to
characterize a scaling function. We will discuss these two
problems.

\vspace*{10pt}
\subsection*{5. Teichm\"uller structures, part I}$\\$

\vspace*{-5pt}

A discussion of the first problem follows a similar idea to that
of Teichm\"uller theory for Riemann surfaces with the help of the
following theorem (refer to~\cite{Ji1,Ji2,Ji6}).

\vspace*{10pt}

\begin{theorem}~\label{qs}
Suppose $f$ and $g$ are two maps in ${\mathcal C}^{1+}_{d}$. Suppose
$h$ is the topological conjugacy between $f$ and $g$. Then $h$ is a
quasisymmetric homeomorphism.
\end{theorem}

A homeomorphism $h$ of $T$ is called quasisymmetric (see~\cite{Al})
if there is a constant $K\geq 1$ such that
$$
K^{-1} \leq \frac{|H(x+t)-H(x)|}{|H(x)-H(x-t)|}\leq K
$$
for all $x\in {\mathbb R}$ and all $t>0$, where $H$ is a lift of $f$
to the real line.

Take $q_{d} (z) =z^{d}$ as a basepoint in ${\mathcal C}^{1+}_{d}$.
For any $f\in {\mathcal C}^{1+}_{d}$, let $h_{f}$ be the conjugacy
from $f$ to $q_{d}$, i.e.,
$$
f\circ h_f= h_f\circ q_{d}.
$$
Thus, we can think of ${\mathcal C}^{1+}_{d}$ as pairs $(f,
h_{f})$. Two pairs satisfy $(f,h_{f}) \sim_{t} (g,h_{g})$ if
$h_{f}\circ h_{g}^{-1}$ is a $C^{1}$-diffeomorphism. Then $\sim_t$
is an equivalence relation. The Teichm\"uller space
$$
{\mathcal T}{\mathcal C}^{1+}_{d}=\{ [(f, h_f)] \;|\; f\in {\mathcal
C}^{1+}_{d}, \; \hbox{with the basepoint $[(q_{d}, id)]$} \}
$$
is the space of all $\sim_{t}$-equivalence classes
$[(f,h_{f})]=[(f,h_{f})]_{t}$ with the basepoint $[(q_{d}, id)]$.
This space has a Teichm\"uller metric $d_{T}(\cdot, \cdot)$ as we
describe below.

We first consider the universal Teichm\"uller space. Let
${\mathcal QS}$ be the set of all quasisymmetric homeomorphisms of
the unit circle $T$ factored by the space of all M\"obius
transformations of the circle. (Then ${\mathcal QS}$ may be
identified with the set of all quasisymmetric homeomorphisms of
the unit circle fixing three points). For any $h\in {\mathcal
QS}$, let ${\mathcal E}_{h}$ be the set of all quasiconformal
extensions of $h$ into the unit disk. Let $K_{\widetilde{h}}$ be
the quasiconformal dilatation of $\widetilde{h}\in {\mathcal
E}_{H}$. Using quasiconformal dilatation, one defines a distance
in ${\mathcal QS}$ by
$$
d_{T}(h_{1}, h_{2}) =\frac{\ 1\ }{\ 2\ }\inf \{ \log
K_{\widetilde{h}_{1}\widetilde{h}_{2}^{-1}}\;|\;
\widetilde{h}_{1}\in {\mathcal E}_{h_{1}},\widetilde{h}_{2} \in
{\mathcal E}_{2}\}.
$$
Here $({\mathcal QS}, d)$ is called the universal Teichm\"uller
space. It is a complete metric space and a complex manifold with
complex structure compatible with the Hilbert transform (see, for
example,~\cite{Al}).

A quasisymmetric homeomorphism $h$ is called symmetric if there is
a bounded positive function $\epsilon (t)$ such that $\epsilon (t)
\rightarrow 0^+$ as $t\rightarrow 0^+$ and
$$
1-\epsilon (t) \leq \frac{|H(x+t)-H(x)|}{|H(x)-H(x-t)|} \leq
1+\epsilon (t)
$$
for any $x$ in ${\mathbb R}$, where $H (x+1)=H(x)+1$ is a lift of
$h$. A $C^{1}$-diffeomorphism of the unit circle is symmetric.
However, a symmetric homeomorphism of the unit circle could be
very singular. Let ${\mathcal S}$ be the subset of ${\mathcal QS}$
consisting of all symmetric homeomorphisms of the unit circle. The
space ${\mathcal S}$ is a closed subgroup of ${\mathcal QS}$. The
topology coming from the metric $d_{T}$ on ${\mathcal QS}$ induces
a topology on the factor space ${\mathcal QS} {\rm \ mod \
}{\mathcal S}$. Given two cosets ${\mathcal S} f$ and ${\mathcal
S} g$ in this factor space, define a metric by
$$
\overline{d}_{T} ({\mathcal S}f, {\mathcal S}g) =\inf_{A, B\in
{\mathcal S}} d(Af, Bg).
$$
The factor space ${\mathcal QS} {\rm \ mod \ }{\mathcal S}$ with
this metric is a complete metric space and a complex manifold. The
topology on $({\mathcal QS} {\rm \ mod \ } {\mathcal S},
\overline{d}_{T})$ is the finest topology which makes the
projection $\pi: {\mathcal QS} \to {\mathcal QS} {\rm \ mod \ }
{\mathcal S}$ continuous, and $\pi$ is also holomorphic. An
equivalent topology can be defined as follows. For any $h\in
{\mathcal QS}$, let $\widetilde{h}$ be a quasiconformal extension
of $h$ to a small neighborhood $U$ of $T$ in the complex plane.
Let
$$
\mu_{\widetilde{h}}(z)
=\frac{\widetilde{h}_{\overline{z}}(z)}{\widetilde{h}_{z}(z)},
\quad k_{\widetilde{h}}=\|\mu_{\widetilde{h}}\|_{\infty} \quad
\hbox{and}\quad  B_{\widetilde{h}}
=\frac{1+k_{\widetilde{h}}}{1-k_{\widetilde{h}}}.
$$
Then the boundary
dilatation is defined as
$$
B_{h} =\inf B_{\widetilde{h}}
$$
where the infimum is taken over all quasiconformal extensions of $h$
near the unit circle. It is known that $h$ is symmetric if and only
if $B_{h}=1$. Define
$$
\widetilde{d} (h_{1},h_{2}) =\frac{1}{2} \log B_{h_{2}^{-1}h_{1}}.
$$
The two metrics $\overline{d}$ and $\widetilde{d}$ on ${\mathcal
QS} {\rm \ mod \ } {\mathcal S}$ are equal. The reader may refer
to~\cite{GS} for this. The Teichm\"uller metric on ${\mathcal
T}{\mathcal C}^{1+}_{d}$ is defined similarly. Let $\tau$ and
$\tau'$ be two points in ${\mathcal T}{\mathcal C}^{1+}_{d}$. Then
$$
d_{T} (\tau, \tau') = \frac{1}{2} \log B_{h_{f}^{-1}\circ h_{g}}
$$
where $(f, h_{f})\in \tau$ and $(g,\tau_{g})\in \tau'$.

Since the space of geometric models can be represented by the
space of scaling functions, {\em the Techm\"uller space can also
be represented by the space of H\"older continuous scaling
functions $S_{f}$ for $f\in {\mathcal C}^{1+}_{d}$ with the
basepoint $S(\kappa)=1/d$.}

\subsection*{6. Characterizations, part I}$\\$

\vspace*{-5pt}

The characterization of a scaling function has been done for $d=2$,
which is the most interesting case for us. First by the definition
of a scaling function, one can easily check that
$$
S(\kappa0)+S(\kappa1)=1, \quad \kappa\in \Sigma^{-}=\Sigma^{-}_{2}.
$$
We call this the summation condition. In addition to this condition,
a scaling function also enjoys another non-trivial condition which
we call the compatibility condition,
$$
\prod_{n=0}^{\infty} \frac{S(\kappa1\overbrace{0\cdots
0}^{n})}{S(\kappa0\underbrace{1\cdots 1}_{n})} =const., \quad
\kappa\in \Sigma^{-}.
$$
Actually, this infinite product converges to the constant
exponentially. (Its general term must tend to $1$ as $n$ goes to
$\infty$. This implies that $S(\cdots 000) =S(\cdots 111)$.) We
showed that the converse is also true as follows.

\vspace*{10pt}

\begin{theorem}~\label{ch1}
Let $S$ be a positive H\"older continuous function on $\Sigma^{-}$.
Then $S$ is the scaling function of a map in ${\mathcal C}^{1+}_{2}$
if and only if $S$ satisfies the summation and compatibility
conditions.
\end{theorem}

\vskip5pt The original proof of this theorem is given
in~\cite{CJQ} and uses the Gibbs measure theory and some
constructions in~\cite{Qu}. A proof without using the Gibbs
measure theory can be founded in~\cite{Ji7}. In the proof, we find
the connection between the scaling function and the solenoid
function and the linear model for a circle endomorphism and use
some constructions in~\cite{CGJ} (refer to Theorem~\ref{ssy}). We
would like to note that the solenoid function and the linear model
are also interesting geometric invariants for a circle
endomorphism. The solenoid function for a circle endomorphism is
defined in~\cite{Su2} and used to describe an affine structure
along leave directions $\Sigma^{+}$ of the solenoid
$\Sigma^{-}\times \Sigma^{+}$. It has been studied in~\cite{PS}.
The linear model for a circle endomorphism is defined
in~\cite{DH}. A linear model can be thought of as a nonlinear
coordinate on the unit circle. In~\cite{DH}, a question about what
kind of nonlinear coordinate can be realized by a smooth expanding
circle endomorphism arose. This question was studied in~\cite{Cu}
by employing some results in quasiconformal theory (refer
to~\cite{Al}). A much simpler understanding was given
in~\cite{Ji7} by employing the naive distortion property
(see~\cite{Ji2}).

Therefore, {\em the Tecihm\"uller space ${\mathcal T}{\mathcal
C}^{1+}_{2}$ is represented by the space of positive H\"older
continuous functions on $\Sigma^{-}$ satisfying the summation and
compatibility conditions with the basepoint $S(\kappa) =1/2$.}

It is clear that the summation condition is true for any $d>2$. This
means that the scaling function $S$ of $f\in {\mathcal C}^{1+}_{d}$
satisfies
$$
S(\kappa0)+S(\kappa1)+\cdots +S(\kappa(d-1)) =1, \quad \kappa\in
\Sigma^{-}=\Sigma^{-}_{d}.
$$
The compatibility condition for $d>2$ should be similar. However,
the proof of the characterization of the scaling function for a
map in ${\mathcal C}^{1+}_{d}$ for $d>2$ should be slightly more
complicated than the case $d=2$, but it is a promising problem.

\vskip10pt
\subsection*{7. Teichm\"uller structures, part II}$\\$

\vspace*{-5pt}

The Teichm\"uller space $({\mathcal T}{\mathcal C}^{1+}_{d},
d_{T}(\cdot, \cdot))$ is not complete. Its completion is an
interesting subject to be studied. A circle endomorphism $f$ of
degree $d$ is called uniformly symmetric if all its inverse branches
for $f^n$, $n=1, 2,\cdots$, are symmetric uniformly. More precisely,
there is a bounded positive function $\epsilon (t)$ with $\epsilon
(t)\to 0^{+}$ as $t\to 0^{+}$ such that
$$
1-\epsilon (t) \leq \frac{|F^{-n} (x+t)
-F^{-n}(x)|}{|F^{-n}(x)-F^{-n}(x-t)|} \leq 1+\epsilon (t), \quad
x\in {\mathbb R},\; t>0,\; n=1, 2, \cdots.
$$
By the naive distortion lemma (see, for example, ~\cite{Ji2}) we
have

\vspace*{10pt}
\begin{proposition} Any map $f\in {\mathcal C}^{1+}$ is uniformly
symmetric.
\end{proposition}

\vskip5pt Let ${\mathcal U}{\mathcal S}$ be the space of all
uniformly symmetric circle endomorphisms of degree $d\geq 2$. The
above proposition says that ${\mathcal C}^{1+}\subset {\mathcal
U}{\mathcal S}$. However, a map in ${\mathcal U}{\mathcal S}$ can
be quite different. For example, it may not be differentiable and
may not be absolutely continuous. However, we have shown that from
the dual point of view, it has a lot of similarity to what we have
studied for a map in ${\mathcal C}^{1+}$. (However, a $C^{1}$
expanding circle endomorphism is very different from what we have
studied for a map in ${\mathcal C}^{1+}$ (see~\cite{Qu}).)

For a fixed $d\geq 2$, let ${\mathcal U}{\mathcal S}_{d}$ be the
space of all uniformly symmetric circle endomorphisms of degree $d$.
For $f\in {\mathcal U}{\mathcal S}_{d}$, it is certainly uniformly
$M$-quasisymmetric for a fixed constant $M>1$, that is,
$$
M^{-1} \leq \frac{|F^{-n} (x+t)
-F^{-n}(x)|}{|F^{-n}(x)-F^{-n}(x-t)|} \leq M, \quad x\in {\mathbb
R},\; t>0,\; n=1, 2, \cdots.
$$
We have:

\vspace*{10pt}
\begin{proposition}
If $f\in {\mathcal U}{\mathcal S}_{d}$, then there is a constant
$C>0$ such that
$$
S(\kappa_{n}) = \frac{|I_{\kappa_{n}}|}{|I_{\sigma^{-}
(\kappa_{n})}|} \geq C
$$
for all finite words $\kappa_{n}=j_{n-1}\cdots j_{1}j_{0}$ of $\{0,
1, \cdots, d-1\}$.
\end{proposition}

The above proposition means that the sequence of nested partitions
$$
\{ \eta_{n}\}_{n=1}^{\infty}=\{ \{ I_{w_{n}}\}_{w_{n}}
\}_{n=1}^{\infty}
$$
has bounded geometry. Using the above proposition and the summation
condition, we have constants $D>0$ and $0<\tau <1$ such that
$$
\varepsilon_{n} \leq D\tau^{n}, \quad n=1, 2, \cdots.
$$
Therefore, just like in the proofs of Theorem~\ref{tc} and
Theorem~\ref{qs}, we have:

\vspace*{10pt}
\begin{theorem}
Any two maps $f, g\in {\mathcal U}{\mathcal S}_{d}$ are
topologically conjugate and, furthermore, the conjugacy $h$ is
quasisymmetric.
\end{theorem}

With this proposition, we can define the Teichm\"uller space for
${\mathcal U}{\mathcal S}_{d}$ as we did for ${\mathcal
C}^{1+}_{d}$. Take $q_{d} (z) =z^{d}$ as a basepoint in ${\mathcal
U}{\mathcal S}_{d}$. For any $f\in {\mathcal U}{\mathcal S}_{d}$,
let $h_{f}$ be the conjugacy from $f$ to $q_{d}$, i.e.,
$$
f\circ h_f= h_f\circ q_{d}.
$$
Thus we can think of ${\mathcal U}{\mathcal S}_{d}$ as pairs $(f,
h_{f})$. Two pairs satisfy $(f,h_{f}) \sim_{t} (g,h_{g})$ if
$h_{f}\circ h_{g}^{-1}$ is symmetric. Then $\sim_t$ is an
equivalence relation. The Teichm\"uller space
$$
{\mathcal T}{\mathcal U}{\mathcal S}_{d}=\{ [(f, h_f)] \;|\; f\in
{\mathcal U}{\mathcal S}_{d}, \; \hbox{with the basepoint $[(q_{d},
id)]$}\}
$$
is the space of all $\sim_{t}$-equivalence classes $[(f, h_f)]=[(f,
h_f)]_{t}$ equipped with a Teichm\"uller metric
$$
d_{T}(\tau, \tau') = \frac{1}{2} \log B_{h_{f}^{-1}\circ h_{g}}
$$
where $(f, h_{f})\in \tau$ and $(g,h_{g})\in \tau'$.

If $f, g\in {\mathcal C}^{1+}_{d}$ and if the conjugacy $h$
between $f$ and $g$ is symmetric, then $h$ must be $C^{1}$. The
reason is that if the conjugacy between $f$ and $g$ is symmetric,
then their scaling functions $S_{f}$ and $S_{g}$ must be the same
(refer to Theorem~\ref{ssy}). Therefore they are $C^{1}$-conjugate
and the conjugacy $h$ must be $C^{1}$ (see Theorem~\ref{sc}). (A
related easy but interesting fact is that the ratio of eigenvalues
$e_{f}(p)$ and $e_{g}(h(p))$ of $f$ and $g$ at corresponding
periodic point $p$ and $h(p)$ determines the local quasisymetric
constant of $h$ at $p$. If $h$ is symmetric, its local
quasisymmetric constant at $p$ is $1$, so the ratio
$e_{f}(p)/e_{g}(h(p))$ is $1$.) This implies that the
Teichm\"uller space ${\mathcal T} {\mathcal C}^{1+}_{d}$ is indeed
a subspace of the Teichm\"uller space ${\mathcal T}{\mathcal
U}{\mathcal S}_{d}$. Furthermore, we have (refer
to~\cite{CJQ,CGJ}):

\vspace*{10pt}
\begin{theorem}
The space $({\mathcal T}{\mathcal U}{\mathcal S}_{d}, d_{T}(\cdot,
\cdot))$ is a complete complex Banach manifold and is the completion
of the space $({\mathcal T} {\mathcal C}^{1+}_{d}, d_{T}(\cdot,
\cdot))$.
\end{theorem}

The local model of the complex Banach manifold can be thought of as
the set of Beltrami coefficients on the upper-half plane ${\mathbb
H}$ (complex $L^\infty$ functions $\mu(z)$ on the upper-half plane
${\mathbb H}$ with $\| \mu (z) \|_{\infty} <1$) such that
$\mu(dz)=\mu(z)$ and $|\mu (z+n)-\mu(z)| \to 0$ uniformly for $n$ as
$\Im(z) \to 0$ (refer to~\cite{Cu,CJQ,CGJ}).

\vspace*{10pt}
\subsection*{8. Geometric models, part II}$\\$

\vspace*{-5pt}

The geometric models of maps in ${\mathcal U}{\mathcal S}_{d}$ can
also be represented by their scaling functions. Two maps $f, g\in
{\mathcal U}{\mathcal S}_{d}$ are called symmetrically conjugate
if the conjugacy between them is symmetric. This is an equivalence
relation which we denote as $f\sim_{sy} g$. The space ${\mathcal
U}{\mathcal S}_{d}/\sim_{sy}$ of geometric models for maps in
${\mathcal U}{\mathcal S}_{d}$ is the space of all equivalence
classes. We have (refer to~\cite{CGJ}):

\vspace*{10pt}

\begin{theorem}~\label{ssy}
Suppose $f\in {\mathcal U}{\mathcal S}_{d}$. Then its scaling
function
$$
S=S_{f} : \Sigma^{-}\to {\mathbb R}^{+}.
$$
exists and is a continuous function. Furthermore, it is a complete
symmetric invariant for ${\mathcal U}{\mathcal S}_{d}$; this means
$f$ and $g$ are symmetrically conjugate if and only if their
scaling functions are the same, i.e., $S_{f}=S_{g}$.
\end{theorem}

Thus ${\mathcal T} {\mathcal U}{\mathcal S}_{d}$ can be
represented by scaling functions $S_{f}$, i.e.,
$$
{\mathcal U}{\mathcal S}_{d}/\sim_{sy} =\{ S_{f}\; |\; f\in
{\mathcal U}{\mathcal S}_{d}\}
$$
and
$$
\big( {\mathcal T} {\mathcal U}{\mathcal S}_{d}=\{ S_{f}\; |\; f\in
{\mathcal U}{\mathcal S}_{d},\; \hbox{with the basepoint
$S=\frac{1}{d}$}\}, d_{T}(\cdot,\cdot) \big).
$$

\vspace*{10pt}

\subsection*{9. Characterizations, part II}$\\$

\vspace*{-5pt}

The characterization of the scaling functions for ${\mathcal
U}{\mathcal S}_{2}$ has been given as

\vspace*{10pt}
\begin{theorem}~\label{ch2}
Let $S$ be a positive continuous function on
$\Sigma^{-}=\Sigma^{-}_{2}$. Then $S$ is the scaling function of a
map in ${\mathcal U}{\mathcal S}_{2}$ if and only if $S$ satisfies
the summation and compatibility conditions.
\end{theorem}

The proof of this theorem can be founded in~\cite{CGJ}. In this
case, the infinite product
$$
\prod_{n=0}^{\infty} \frac{S(\kappa1\overbrace{0\cdots
0}^{n})}{S(\kappa0\underbrace{1\cdots 1}_{n})} =const., \quad
\kappa\in \Sigma^{-},
$$
in the compatibility condition converges uniformly to a constant.

Therefore, {\em the Tecihm\"uller space ${\mathcal T}{\mathcal
U}{\mathcal S}_{2}$ is represented by the space of positive
continuous functions on $\Sigma^{-}$ satisfying the summation and
compatibility conditions.}

Just as in the end of \S6, the characterization of a scaling
function of a map in ${\mathcal U}{\mathcal S}_{d}$ should be
slightly more complicated than the case $d=2$, but it is a
promising problem.

\vskip10pt
\subsection*{10. Invariant measures and dual invariant measures}$\\$

\vspace*{-5pt}

Consider the symbolic dynamical system $(\Sigma^{+}, \sigma^{+})$
and a positive H\"older continuous function $\psi=\psi (w)$. The
standard Gibbs theory (refer to~\cite{Bo} or~\cite{FJ}) implies that
there is a number $P=P(\log \psi)$ called the pressure and a
$\sigma^{+}$-invariant probability measure $\mu_{+}= \mu_{+,\psi}$
such that
$$
C^{-1} \leq \frac{\mu_{+} ([i_{0}\cdots i_{n-1}])}{\exp (-Pn +
\sum_{i=0}^{n-1} \log \psi ((\sigma^{+})^{i} (w)))}\leq C
$$
for any left cylinder $[i_{0}\cdots i_{n-1}]$ and any $w=i_{0}\cdots
i_{n-1}\cdots \in [i_{0}\cdots i_{n-1}]$, where $C$ is a fixed
constant. Here $\mu_{+}$ is a $\sigma^{+}$-invariant measure means
that
$$
\mu_{+} ((\sigma^{+})^{-1} (A)) =\mu_{+} (A)
$$
for all Borel sets of $\Sigma^{+}$. A $\sigma^{+}$-invariant
probability measure satisfying the above inequalities is called the
Gibbs measure with respect to the given potential function $\log
\psi$.

Two positive H\"older continuous functions $\psi_{1}$ and $\psi_{2}$
are said to be cohomologous equivalent if there is a continuous
function $u=u(w)$ on $\Sigma^{+}$ such that
$$
\log \psi_{1}(w) -\log \psi_{2} (w) = u (\sigma^{+} (w)) -u(w).
$$
If two functions are cohomologous to each other, they have the same
Gibbs measure. Therefore, the Gibbs measure can be thought of as a
representation of a cohomologous class.

The Gibbs measure is also an equilibrium state. Consider the
measure-theoretical entropy $h_{\mu_{+}}(\sigma^{+})$. Since the
Borel $\sigma$-algebra of $\Sigma^{+}$ is generated by all left
cylinders, then $h_{\mu_{+}} (\sigma^{+}) $ can be calculated as
$$
h_{\mu_{+}} (\sigma^{+})  =\lim_{n\to \infty} {1\over n}
\sum_{w_{n}} \Big( -\mu_{+} ([w_{n}]) \log \mu_{+} ([w_{n}])\Big)
$$
$$
= \lim_{n\to \infty} \sum_{w_{n}} \Big( -\mu_{+} ([w_{n}]) \log
\big(\frac{\mu_{+}(
[w_{n}])}{\mu_{+}(\sigma^{+}([w_{n}]))}\big)\Big),
$$
where $w_{n}$ runs over all words $w_{n}=i_{0}\cdots i_{n-1}$ of $\{
0, 1,\cdots, d-1\}$ of length $n$. Then $\mu_{+}$ is an equilibrium
state in the sense that
$$
P(\log \psi) = h_{\mu_{+}} (\sigma^{+}) + \int_{\Sigma^{+}} \log
\psi (w) d\mu (w) =\sup \{ h_{\nu}(\sigma^{+}) +\int_{\Sigma^{+}}
\log \psi (w) d\nu(w)\}
$$
where $\nu$ runs over all $\sigma^{+}$-invariant probability
measures. The measure $\mu_{+}$ is unique in this case.

There is a natural way to transfer a $\sigma^{+}$-invariant
probability measure $\mu_{+}$ (not necessarily a Gibbs measure) to a
$\sigma^{-}$-invariant probability measure $\mu_{-}$ as follows.
Given any right cylinder $[j_{n-1}\cdots j_{0}]_{r}$ in
$\Sigma^{-}$, let $i_{0}\cdots i_{n-1}=j_{n-1}\cdots j_{0}$ define a
left cylinder
$$
[i_{0}\cdots i_{n-1}]_{l} =\{ w'=i_{0}'\cdots i_{n-1}'i_{n}'\cdots
\; |\; i_{0}'=i_{0}, \cdots, i_{n-1}'=i_{n-1}\}.
$$
Then define
$$
\mu_{-} ([j_{n-1}\cdots j_{0}]_{r})= \mu_{+}([i_{0}\cdots
i_{n-1}]_{l}).
$$
Then
$$
\mu_{-} ([j_{n-1}\cdots j_{0}]_{r})= \mu_{+}([i_{0}\cdots
i_{n-1}]_{l})= \mu_{+}((\sigma^{+})^{-1}([i_{0}\cdots i_{n-1}]_{l}))
$$
$$
=\mu_{+}(\cup_{i=0}^{d-1} [ii_{0}\cdots i_{n-1}]_{l}) =
\sum_{i=0}^{d-1} \mu_{+}([ii_{0}\cdots i_{n-1}]_{l})
=\sum_{j=0}^{d-1} \mu_{-}([jj_{n-1}\cdots j_{0}]).
$$
This implies that $\mu_{-}$ satisfies the finite additive law for
all cylinders, i.e., if $A_{1}$, $\cdots$, $A_{k}$ are finitely many
pairwise disjoint right cylinders in $\Sigma^{-}$, then
$$
\mu_{-} (\cup_{l=1}^{k} A_{k}) =\sum_{l=1}^{k} \mu_{-} (A_{l}).
$$
Also $\mu_{-}$ satisfies the continuity law in the sense that if
$\{ A_{n}\}_{n=1}^{\infty}$ is a decreasing sequence of cylinders
and tends to the empty set (this means $A_{n+1}\subset A_{n}$ and
$\cap_{n=1}^{\infty}A_{n}=\emptyset$), then $\mu_{-} (A_{n})$
tends to zero as $n$ goes to $\infty$. The reason is that since a
cylinder of $\Sigma^{-}$ is a compact set, a decreasing sequence
of cylinders tending to the empty set must be eventually all
empty. The Borel $\sigma$-algebra in $\Sigma^{-}$ is generated by
all right cylinders. So $\mu_{-}$ extends to measure on
$\Sigma^{-}$. We have the following proposition.

\vspace*{10pt}
\begin{proposition}
$\mu_{-}$ is a $\sigma^{-}$-invariant probability measure.
\end{proposition}

\begin{proof} We have seen that $\mu^{-}$ is a measure on
$\Sigma^{-} $. Since $\mu_{-}(\Sigma^{-})=1$, it is a probability
measure. For any right cylinder $[j_{n-1}\cdots j_{0}]_{r}$,
$$
\mu_{-} ((\sigma^{-})^{-1} ([j_{n-1}\cdots j_{0}]_{r}) = \mu_{-}
(\cup_{j=0}^{d-1} [j_{n-1}\cdots j_{0}j]_{r})
$$
$$
=\sum_{j=0}^{d-1} \mu_{-} ([j_{n-1}\cdots j_{0}j]_{r}) =
\sum_{i=0}^{d-1} \mu_{+} ([i_{0}\cdots i_{n-1}i]_{l})
$$
$$
= \mu_{+} (\cup_{i=0}^{d-1} [i_{0}\cdots i_{n-1}i]_{l}) =
\mu_{+}([i_{0}\cdots i_{n-1}]_{l}) =\mu_{-}([j_{n-1}\cdots
j_{0}]_{r}).
$$
So $\mu_{-}$ is $\sigma^{-}$-invariant.
\end{proof}

We call $\mu_{-}$ a dual invariant measure. A natural question now
is as follows. {\em Is a dual invariant measure a Gibbs measure with
respect to some continuous or H\"older continuous function on
$\Sigma^{-}$?}

A more interesting geometric question is the following. Consider a
metric induced from the dual invariant measure $\mu_{-}$ (in the
case that $\mu_{-}$ is supported on the whole $\Sigma^{-}$ and is
non-atomic), that is,
$$
d(\kappa, \kappa') =\mu_{-}([j_{n-1}\cdots j_{0}])
$$
where $[j_{n-1}\cdots j_{0}]$ is the smallest right cylinder
containing both $\kappa=\cdots j_{n}j_{n-1}\cdots j_{0}$ and
$\kappa'=\cdots j_{n}' j_{n-1}\cdots j_{0}$, $j_{n}\neq j_{n}'$.
{\em Is $\sigma^{-}$ differentiable with a continuous or H\"older
continuous derivative under this metric?} More precisely, {\em does
the limit
$$
\frac{d\sigma^{-}}{dx} (\kappa)= \lim_{n\to \infty}
\frac{\mu_{-}(\sigma^{-}([j_{n-1}\cdots
j_{1}j_{0}]))}{\mu_{-}([j_{n-1}\cdots j_{1}j_{0}])} = \lim_{n\to
\infty} \frac{\mu_{-}([j_{n-1}\cdots j_{1}])}{\mu_{-}([j_{n-1}\cdots
j_{1} j_{0}])}
$$
exist for every $\kappa=\cdots j_{n-1}\cdots j_{1}j_{0} \in
\Sigma^{-}$? If it exists, is the limiting function continuous or
H\"older continuous on $\Sigma^{-}$?}

Actually, there is a measure-theoretical version related to these
questions. I will first give a brief review of this theory.

\vskip10pt
\subsection*{11. $g$-measures}$\\$

\vspace*{-5pt}

Let $X$ be $\Sigma^{-}$ (or $\Sigma^{+}$) and let $f$ be
$\sigma^{-}$ (or $\sigma^{+}$). Let ${\mathcal B}$ be the Borel
$\sigma$-algebra of $X$. Let ${\mathcal M}(X)$ be the space of all
finite Borel measures on $X$. Let ${\mathcal M}(X, f)$ be the space
of all $f$-invariant probability measures in ${\mathcal M}(X)$. Let
${\mathcal C} (X)$ be the space of all continuous real functions on
$X$. Then ${\mathcal M}(X)$ is the dual space of ${\mathcal C}(X)$.
Denote
$$
<\phi, \mu> =\int_{X} \phi (x) d\mu, \quad \hbox{$\phi \in {\mathcal
C}(X)$ and $\mu\in {\mathcal M}(X)$}.
$$

A real non-negative continuous function $\psi$ on $X$ is called a
$g$-function (the historic reason to call such a function a
$g$-function is because of Keane's paper~\cite{Ke}) if
$$
\sum_{fy=x} \psi(y) =1.
$$

For a function $\psi$, define the transfer operator ${\mathcal
L}_{\psi}$ from ${\mathcal C}(X)$ into itself as
$$
{\mathcal L}_{\psi} \phi (x) =\sum_{f(y)=x} \phi (y) \psi(y), \quad
\phi \in {\mathcal C}(X).
$$
One can check that ${\mathcal L}_{\psi} \phi = {\mathcal L}_{1}
(\psi\phi)$ and if $\psi$ is a $g$-function, then ${\mathcal
L}_{\psi}1=1$. Let ${\mathcal L}_{\psi}^{*} $ be the dual operator
of ${\mathcal L}_{\psi}$, that is, ${\mathcal L}_{\psi}^{*}$ is the
operator from ${\mathcal M}(X)$ into itself satisfying
$$
<\phi, {\mathcal L}_{\psi}^{*}\mu> =<{\mathcal L}_{\psi}\phi, \mu>,
\quad \forall \phi \in {\mathcal C}(X) \; \hbox{and}\; \forall \mu
\in {\mathcal M}(X).
$$

Suppose $\psi$ is a $g$-function. Then a probability measure $\mu\in
{\mathcal M}(X)$ is called a $g$-measure if it is a fixed point of
${\mathcal L}_{\psi}$, that is,
$$
{\mathcal L}_{\psi}^{*} \mu =\mu.
$$
A $g$-measure is a $f$-invariant measure because
$$
\mu (f^{-1}(B)) = <1_{f^{-1}(B)}, \mu>= < 1_{B}\circ f, {\mathcal
L}_{\psi}^{*} \mu>
$$
$$
= <{\mathcal L}_{\psi} 1_{B}\circ f, \mu> =<1_{B}, \mu>=\mu (B),
\; \forall B\in {\mathcal B}.
$$

For any $\mu\in {\mathcal M}(X)$, let $\widetilde{\mu} ={\mathcal
L}_{1}^{*}\mu$.

\vspace*{10pt}
\begin{proposition}
$$
\widetilde{\mu} (B) = \sum_{j=0}^{d-1} \mu (f(B\cap [j]))
$$
where $B$ is any Borel subset in ${\mathcal B}$ and $[j]$ is the
right cylinder of $j$. Moreover, if $\mu\in {\mathcal M}(X,f)$,
$\mu$ is absolutely continuous with respect to $\widetilde{\mu}$.
\end{proposition}

\begin{proof}
For any Borel subset $B\in {\mathcal B}$,
$$
\widetilde{\mu} (B) = <1_{B}, {\mathcal L}_{1}^{*}\mu> =
<{\mathcal L}_{1} 1_{B}, \mu>.
$$
But
$$
{\mathcal L}_{1} 1_{B} (x) = \sum_{j=0}^{d-1} 1_{B}(xj) =
\sum_{j=0}^{d-1} 1_{f(B\cap [j])} (x).
$$
So we have that
$$
\widetilde{\mu} (B) = \sum_{j=0}^{d-1} \mu (f(B\cap [j])).
$$

If $\mu$ is $f$-invariant, then we have
$$
\widetilde{\mu} (B) = \sum_{j=0}^{d-1} \mu (f(B\cap [j])) =
\sum_{j=0}^{d-1} \mu (f^{-1}(f(B\cap [j])))\geq \sum_{j=0}^{d-1}
\mu (B\cap [j]) =\mu (B).
$$
Therefore, $\mu(B)=0$ whenever $\widetilde{\mu}(B) =0$. So $\mu$
is absolutely continuous with respect to $\widetilde{\mu}$.
\end{proof}

Suppose $\mu\in {\mathcal M}(X,f)$. Then $\mu$ is absolutely
continuous with respect to $\widetilde{\mu}$. So the
Radon-Nikod\'ym derivative
$$
D_{\mu}(x)=\frac{d\mu} {d\widetilde{\mu}} (x), \quad
\widetilde{\mu}-a.e.\; x
$$
of $\mu$ with respect to $\widetilde{\mu}$ exists
$\widetilde{\mu}$-a.e. and is a $\widetilde{\mu}$-measurable
function. We would like to note that $\widetilde{\mu}$ may not be
absolutely continuous with respect to $\mu$.

The following theorem was proved by Leddraper in~\cite{Le} and was
used by Walters in~\cite{Wa} in the study of a generalized version
of Ruelle's theorem.

\vspace*{10pt}
\begin{theorem}
Suppose $\psi$ is a $g$-function and $\mu\in {\mathcal M}(X)$ is a
probability measure. The followings are equivalent:
\begin{itemize}
\item[i)] $\mu$ is a $g$-measure, i.e., ${\mathcal L}_{\psi}^{*}\mu
=\mu$.
\item[ii)] $\mu\in {\mathcal M}(X,f)$ and $D_{\mu} (x) =\psi(x)$ for
$\widetilde{\mu}$-a.e. $x$.
\item[iii)] $\mu\in {\mathcal M}(X,f)$ and
$$
E[\phi | f^{-1}({\mathcal B})] (x) = {\mathcal L}_{\psi} \phi
(fx)=\sum_{fy=fx} \psi (y) \phi(y), \; \hbox{for $\mu$-a.e.\; x}
$$
where $E[\phi|f^{-1}({\mathcal B})]$ is the conditional expectation
of $\phi$ with respect to $f^{-1}({\mathcal B})$.
\item[iv)] $\mu\in {\mathcal M}(X, f)$ and is an equilibrium state
in the meaning that
$$
0=h_{\mu}(f) +\int_{X} \log \psi\; d\mu = \sup \{ h_{\nu}(f) +
\int_{X} \log \psi\; d\nu \; |\; \nu \in {\mathcal M}(X,f)\}.
$$
(Note that the pressure $P(\log \psi)=0$ for a $g$-function $\psi$.)
\end{itemize}
\end{theorem}

For any $\sigma^{+}$-invariant probability measure $\mu_{+}$, let
$\mu_{-}$ be the dual $\sigma^{-}$-invariant probability measure
which we have constructed in the previous section. Then we have a
$\widetilde{\mu}_{+}$-measurable function
$$
D_{\mu_{+}} (w) =\lim_{n\to \infty} \frac{\mu_{+}
([i_{0}i_{1}\cdots i_{n-1}])}{\mu_{+}([i_{1}\cdots i_{n-1}])},
\quad \hbox{for $\widetilde{\mu}_{+}$-a.e. $w=i_{0}i_{1}\cdots
i_{n-1}\cdots $}
$$
and a $\widetilde{\mu}_{-}$-measurable function
$$
D_{\mu_{-}} (\kappa) =\lim_{n\to \infty} \frac{\mu_{-}
([j_{n-1}\cdots j_{1}j_{0}])}{\mu_{-}([j_{n-1}\cdots j_{1}])},
\quad \hbox{for $\widetilde{\mu}_{-}$-a.e. $\kappa=\cdots
j_{n-1}\cdots j_{0}$}.
$$
Now the question related to those at the end of the previous section
is as follows. Can we extend $D_{\mu_{-}}$ as well as $D_{\mu_{+}}$
to a continuous $g$-function or a H\"older continuous $g$-function?

The Borel $\sigma$-algebra of $\Sigma^{+}$ (or of $\Sigma^{-}$) is
generated by all left cylinders (or all right cylinders). The
measure-theoretical entropy $h_{\mu_{+}}(\sigma^{+})$ can be
calculated as
$$
h_{\mu_{+}} (\sigma^{+})  =\lim_{n\to \infty} {1\over n}
\sum_{w_{n}}\Big( - \mu_{+} ([w_{n}]) \log \mu_{+} ([w_{n}])\Big)
$$
$$
= \lim_{n\to \infty} \sum_{w_{n}} \Big( -\mu_{-} ([w_{n}]) \log
\big(\frac{\mu_{-}(
[w_{n}])}{\mu_{-}(\sigma^{+}([w_{n}]))}\big)\Big),
$$
where $w_{n}$ runs over all words $w_{n}=i_{0}\cdots i_{n-1}$ of $\{
0, 1,\cdots, d-1\}$ of length $n$. The measure-theoretical entropy
$h_{\mu_{-}}(\sigma^{-})$ can be calculated as
$$
h_{\mu_{-}} (\sigma^{-})  =\lim_{n\to \infty} {1\over n}
\sum_{\kappa_{n}}\Big( - \mu_{-} ([\kappa_{n}]) \log \mu_{-}
([\kappa_{n}])\Big)
$$
$$
= \lim_{n\to \infty} \sum_{\kappa_{n}} \Big( -\mu_{-} ([\kappa_{n}])
\log \big(\frac{\mu_{-}(
[\kappa_{n}])}{\mu_{-}(\sigma^{-}([\kappa_{n}]))}\big)\Big),
$$
where $\kappa_{n}$ runs over all words $\kappa_{n}=j_{n-1}\cdots
j_{0}$ of $\{ 0, 1,\cdots, d-1\}$ of length $n$. We would like to
know when is $\mu_{+}$ (or $\mu_{-}$) an equilibrium state? We have
studied these questions for ${\mathcal C}^{1+}$ and for ${\mathcal
U}{\mathcal S}$.

\subsection*{12. Geometric Gibbs measures and dual geometric Gibbs
measures} $\\$

\vspace*{-5pt}

Consider $f\in {\mathcal C}^{1+}$. Then $1/f'(x)$ can be lifted to
a positive H\"older continuous function $\psi (w) =\psi_f (w)
=1/f'(\pi_{+} (w))$ on the symbolic space $\Sigma^{+}$. By
thinking of $\log \psi $ as a potential function for the dynamical
system $(\Sigma^{+}, \sigma^{+})$, there is a unique
$\sigma^{+}$-invariant measure $\mu_{+}= \mu_{+,\psi}$ as we have
mentioned in the previous section such that
$$
C^{-1} \leq \frac{\mu_{+} ([i_{0}\cdots
i_{n-1}])}{\prod_{i=0}^{n-1}\psi ((\sigma^{+})^{i} (w))}\leq C
$$
for any left cylinder $[i_{0}\cdots i_{n-1}]$ and any $w=i_{0}\cdots
i_{n-1}\cdots \in [i_{0}\cdots i_{n-1}]$, where $C$ is a fixed
constant. (Note that $P=P(\log \psi)=0$ in this case.)

The geometric model $[f]_s$ in ${\mathcal C}^{1+}$ can also be
represented by the Gibbs measure $\mu_{+}$ with respect to $\psi (w)
= 1/ f'(\pi_{+} (w))$. The reason is that any $g\in [f]_{s}$ is
smoothly conjugate to $f$, so there is a $C^{1}$ diffeomorphism $h$
of $T$ such that $f(h(x)) =h(g (x))$. Then $f'(h(x)) h'(x) =
h'(g(x)) g'(x)$. Therefore,
$$
\log \psi_{f}(w) -\log \psi_{g} (w) = \log h' (w) - \log
h'(\sigma^{+} (w)).
$$
So $\psi_{g}$ and $\psi_{f}$ are cohomologous to each other. We
call this $\mu_{+}$ a geometric Gibbs measure because it enjoys
the following geometric property too: The push-forward measure
$\mu =(\pi_{+})_{*}\mu_{+}$ is a smooth $f$-invariant measure.
This means there is a continuous function $\rho$ on $T$ such that
$$
\mu (A) =\int_{A} \rho (x) dx, \quad \hbox{for all Borel subsets $A$
on $T$}.
$$

There is another way to find the density $\rho$. First it is a
standard method to find an invariant measure for a dynamical system
$f$. Let $\nu_{0}$ be the Lebesgue measure. Consider the
push-forward measure $\nu_{n}= (f^{n})_{*}\nu_{0}$ by the $n^{th}$
iterates of $f$. Sum up these measures to get
$$
\mu_{n} = \frac{1}{n} \sum_{k=0}^{n-1} \nu_{n}
$$
Any limit $\mu$ of a subsequence of $\{ \mu_{n}\}$ will be an
$f$-invariant measure. Since we start with an $f\in {\mathcal
C}^{1+}$, we can prove that the sequence $\{\mu_{n}\}$ is actually
convergent in $C^{1}$ topology. This means that each $\nu_{n}=
(f^{n})_{*}\nu_{0}$ has a H\"older continuous density
$$
\rho_{n} (x) =\sum_{f^{n}(y)=x} \frac{1}{(f^{n})'(y)}.
$$
Following the theory of transfer operators (refer
to~\cite{Ji8,FJ}), $\rho_{n}(x)$ converges uniformly to a
continuous function $\rho (x)$. The density of $\mu_{n}$ is just
$$
\frac{1}{n} \sum_{k=0}^{n-1} \rho_{n}.
$$
So it also converges to $\rho$ uniformly. Thus $\mu (A) =\int_{A}
\rho (x) dx$ is the limit of $\mu_{n}$ and is a smooth $f$-invariant
measure.

Let $y= h(x) =\mu([0,x])$. Then $y=h(x)$ is a $C^{1}$-diffeomorphism
of $T$. Let
$$
g(y) = h\circ f\circ h^{-1} (y), \quad x= h^{-1}(y)
$$
(Note that $g$ here means a circle endomorphism not a $g$-function!)
Then $g$ preserves the Lebesgue measure $dy$ (which means that
$g_{*}(dy)=dy$, or equivalently, the Lebesgue measure is
$g$-invariant). Since the Lebesgue measure is an ergodic
$g$-invariant measure, $g$ is unique in the geometric model
$[f]_{s}.$

By considering $\psi(w) =1/g'(\pi (w))$, it is a $g$-function on
$\Sigma^{+}$ and $\mu_{+}$ is a $g$-measure. Thus $\mu_{+}$ is an
equilibrium state. It follows that $\mu$ is also an equilibrium
state, that is,
$$
0=P (-\log f'(x))= h_{\mu} (f) -\int_{T} \log f'(x) d\mu = h_{\mu}
(f) -\int_{T} \log f'(x) \rho (x) dx
$$
$$
=\sup \{ h_{\nu}(f) -\int_{T} \log f'(x) d\nu \; |\; \hbox{$\nu$ is
an $f$-invariant propbability measure}\}
$$
$$
= h_{Leb} (g) -\int_{T} \log g'(y) dy,
$$
where $h_{Leb}(g)$ denotes the measure-theoretical entropy with
respect to the Lebesgue measure. The equilibrium state $\mu$ is
unique in this case.

Now by considering the dual invariant measure $\mu_{-}$ for this
geometric Gibbs measure $\mu_{+}$, we have:

\vspace*{10pt}
\begin{theorem}~\label{der1}
Suppose $f\in {\mathcal C}^{1+}$. Consider $\Sigma^{-}$ with the
metric $d(\cdot, \cdot)$ induced from $\mu_{-}$. Then the right
shift $\sigma^{-}$ is a $C^{1+}$ differentiable with respect to $d$.
The derivative is one over the scaling function $S_{f}$, i.e.,
$$
\frac{d \sigma^{-}}{dx} (\kappa) = \frac{1}{S(\kappa)}, \quad
\kappa\in \Sigma^{-}.
$$
\end{theorem}

The proof follows the proof of Theorem~\ref{sc} and the definition
of $\mu_{-}$. Note that by the definition of the derivative for
$\kappa=\cdots j_{n-1}\cdots j_{1}j_{0}\in \Sigma^{-}$,
$$
\frac{d \sigma^{-}}{dx} (\kappa)= \lim_{n\to \infty}
\frac{\mu_{-}(\sigma^{-} ([j_{n-1}\cdots
j_{1}j_{0}]))}{\mu([j_{n-1}\cdots j_{1}j_{0}])} =\lim_{n\to
\infty}\frac{\mu_{-}([j_{n-1}\cdots j_{1}])}{\mu_{-}([j_{n-1}\cdots
j_{1}j_{0}])}.
$$
The theorem says that it equals to $1/S(\kappa)$ pointwise.
Moreover, this convergence is exponentially fast. Then, following
the fact that $\Sigma^{-}$ is a compact space, we have
automatically the Gibbs inequality that
$$
C^{-1} \leq \frac{\mu_{-} ([j_{n-1}\cdots j_{0}])}{\prod_{l=0}^{n-1}
S((\sigma^{-})^{l} (\kappa))} \leq C
$$
for any right cylinder $[j_{n-1}\cdots j_{0}]$ and any $\kappa$ in
this cylinder, where $C>0$ is a fixed constant. Thus $\mu_{-}$ is a
Gibbs measure with respect to the potential function $\log S_{f}$.
We call $\mu_{-}$ a dual geometric Gibbs measure.

\vspace*{10pt}

\begin{corollary}
The dual geometric Gibbs measure $\mu_{-}$ is a $g$-measure with
respect to the $g$-function $S_{f}$ whose pressure $P(\log S_{f})
=0$. Moreover, $D_{\mu_{-}} =S_{f}$ for $\widetilde{\mu}$-a.e.
$\kappa$ and $\mu_{-}$ is a unique equilibrium state in the sense
that
$$
0= P (\log S)= h_{\mu_{-}} (\sigma^{-}) +\int_{\Sigma^{-}} \log
S(\kappa) d\mu_{-}(\kappa)
$$
$$
=\sup \big\{ h_{\nu } (\sigma^{-}) +\int_{\Sigma^{-}} \log S(\kappa)
d\nu (\kappa)\; |\; \hbox{$\nu$ is a $\sigma^{-}$-invariant
measure}\big\}.
$$
\end{corollary}

So following Theorem~\ref{ch1} and Theorem~\ref{der1}, we have:

\vspace*{10pt}

\begin{theorem}~\label{gder1}
Suppose $S (\kappa)$ is a H\"older continuous function on
$\Sigma^{-}$ satisfying the summation and compatibility conditions.
Then there is a unique measure $\mu_{-}$ and the metric $d(\cdot,
\cdot)$ induced from $\mu_{-}$ on $\Sigma^{-}$ such that $1/S$ is
the derivative of the right shift $\sigma^{-}$ with respect to this
metric. Moreover, $\mu_{-}$ is an equilibrium state for the
dynamical system $\sigma^{-}: (\Sigma^{-}, d(\cdot, \cdot)) \to
(\Sigma^{-}, d(\cdot, \cdot))$ from a metric space into itself with
the potential $\log S$.
\end{theorem}

\subsection*{13. Dual geometric Gibbs type measures}$\\$
\vspace*{-5pt}

A map $f\in {\mathcal U}{\mathcal S}$ may not be differentiable
everywhere (it may not be even be absolutely continuous). There is
no suitable Gibbs theory to be used in the study of geometric
properties of a $\sigma^{+}$-invariant measure. We have turned to
the dual symbolic dynamical system $(\Sigma^{-}, \sigma^{-})$ and
produced a similar dual geometric Gibbs type measure theory.

An $f$-invariant measure $\mu$ can be found as we did in the
previous section. Let $\nu_{0}$ be the Lebesgue measure. Consider
the push-forward measure $\nu_{n}=(f^{n})_{*}\nu_{0}$ and sum them
up to get
$$
\mu_{n} =\frac{1}{n} \sum_{i=0}^{n-1} \nu_{n}.
$$
Take a weak limit $\mu$ of a subsequence of $\{\mu_{n}\}$. Then
$\mu$ is an $f$-invariant measure.

Each $h_{n}(x)=\mu_{n}([0,x])$ defines a homeomorphism on $T$.
Since $f$ is uniformly symmetric, the sequence $\{ h_{n}\}$ is
also uniformly symmetric. The space of all quasisymmetric
homeomorphisms with a fixed quasisymmetric constant is a normal
family (refer to~\cite{Al}). So there is a subsequence of $\{
h_{n}\}$ that converges uniformly to a function which is a
symmetric homeomorphism $h(x)$ in this case. Furthermore, we have
$$
h(x) =\mu([0,x]).
$$
Moreover, by considering
$$
g=h\circ f\circ h^{-1},
$$
we see that $g$ is a uniformly symmetric circle endomorphism in
the geometric model $[f]_{sy}$ preserving the Lebesgue measure.

We can lift $\mu$ to $\Sigma^{+}$ to get a $\sigma^{+}$-invariant
measure $\mu_{+}$ as follows. For any finite word $w_{n}=i_{0}\cdots
i_{n-1}$, consider the cylinder $[w_{n}]$. Define
$$
\mu_{+} ([w_{n}]) =\mu (I_{w_{n}}),
$$
where $I_{w_{n}}$ is the interval in the interval system labeled by
$w_{n}$. One can check that it satisfies the finite additive law and
the continuity law. So it can be extended to a
$\sigma^{+}$-invariant probability measure $\mu_{+}$ on $\Sigma^{+}$
such that $(\pi_{+})_{*} \mu_{+}=\mu$. For $\mu_{+}$, we can
construct its dual invariant measure $\mu_{-}$ on $\Sigma^{-}$ as we
did in the previous two sections. Then we have the following
geometric Gibbs type property as we had before in the smooth case:

\vspace*{10pt}

\begin{theorem}~\label{der2}
Suppose $f\in {\mathcal U}{\mathcal S}$. Consider $\Sigma^{-}$ with
the metric $d(\cdot, \cdot)$ induced from $\mu_{-}$. Then the right
shift $\sigma^{-}$ is $C^{1}$ differentiable. The derivative is one
over the scaling function $S_{f}$, i.e.,
$$
\frac{d \sigma^{-}}{dx} (\kappa) = \frac{1}{S_{f}(\kappa)}, \quad
\kappa\in \Sigma^{-}.
$$
\end{theorem}

The proof of the theorem follows the proof of Theorem~\ref{ssy}
and the definition of $\mu_{-}$.

\vspace*{10pt}
\begin{definition}
Suppose $\psi(\kappa)$ is a positive continuous function on
$\Sigma^{-}$. A $\sigma^{-}$-invariant measure $\nu$ is called a
geometric Gibbs type measure with the potential $-\log
\psi(\kappa)$ if
$$
\lim_{n\to \infty} \frac{\nu ([j_{n-1}\cdots
j_{1}j_{0}])}{\nu([j_{n-1}\cdots j_{1}])} =\psi (\kappa), \quad
\forall\; \kappa=\cdots j_{n-1}\cdots j_{0}\in \Sigma^{-}.
$$
\end{definition}

\vspace*{10pt}

\begin{corollary}
The measure $\mu_{-}$ in Theorem~\ref{der2} is a geometric Gibbs
type measure with the potential $\log S_{f} $. Furthermore,
$\mu_{-}$ is a $g$-measure with respect to the $g$-function
$S_{f}$ and $D_{\mu_{-}}(\kappa) =S_{f}(\kappa)$ for
$\widetilde{\mu}_{-}$-a.e. $\kappa$. Moreover, $\mu_{-}$ is an
equilibrium state in the sense that
$$
0= P (-\log S)= h_{\mu_{-}} (\sigma^{-}) +\int_{\Sigma^{-}} \log
S(\kappa) d\mu_{-} (\kappa)
$$
$$
=\sup \{ h_{\nu} (\sigma^{-}) +\int_{\Sigma^{-}} \log S(\kappa)
d\nu (\kappa)\}
$$
where $\nu$ runs over all $\sigma^{-}$-invariant probability
measures.
\end{corollary}

So following Theorem~\ref{ch2} and Theorem~\ref{der2}, we have:

\vspace*{10pt}
\begin{theorem}~\label{gder2}
Suppose $S(\kappa)$ is a continuous function on $\Sigma^{-}$
satisfying the summation and compatibility conditions. Then there
is a geometric Gibbs type measure $\mu_{-}$ with the potential
$\log S$. Moreover, $\mu_{-}$ is an equilibrium state for the
dynamical system $\sigma^{-}$ with the potential $\log S$.
\end{theorem}

The Teichm\"uller space ${\mathcal T}{\mathcal U}{\mathcal S}$ is
a complex Banach manifold and consists of certain positive
functions on the dual symbolic space $\Sigma^{-}$. It is an
interesting problem now to study the change of $\mu_{-}$ when the
potential $\log S$ is changed in the manifold. The reader who is
interested in this direction may refer to~\cite{JR,Ru} for some
results about differentiating the absolutely continuous invariant
measure of a map with respect to this map.

\vskip20pt

\noindent {\em Acknowledgement:} During my research in this
direction and related topics, I have had many interesting
conversations with Professors Guizhen Cui, Aihua Fan, Fred
Gardiner, Jihua Ma, Sudeb Mitra, Huyi Hu, Anthony Quas, David
Ruelle, and Sheldon Newhouse. Professor Dennis Sullivan introduced
me to this interesting topic. I would like to express my thanks to
them all.

\newpage

\bibliographystyle{amsalpha}

\end{document}